\documentclass[12pt]{article}
\usepackage{amsmath}\usepackage{amssymb}\usepackage{amscd}

\newcounter{theorem}\makeatletter
\@addtoreset{theorem}{section}\makeatother

\newtheorem{prop}[theorem]{Proposition}
\newtheorem{lem}[theorem]{Lemma}\newtheorem{def-lem}[theorem]{Definition-Lemma}
\newtheorem{def-prop}[theorem]{Definition-Proposition}

\newtheorem{cor}[theorem]{Corollary}

\newcommand{\rf}[1]{(\ref{#1})}

\def\A{\mathcal{A}}
\def\Ab{\bar{\A}}
\def\Ac{\tilde{\A}}

\def\ba{\begin{eqnarray}}
\def\be{\begin{equation}}
\def\C{\mathbb{C}}

\def\der{\bar{\partial}}
\def\Dif{\text{Diff}}

\def\Diff{\text{Diff}_{\h}}

\def\DR{\mathcal{D}\!\,(\gl)}
\def\DRn{\mathcal{D}\!\,(\gln)}
\newcommand{\e}[2]{e_{#1 #2}}
\newcommand{\E}[2]{e_{#1 #2}}
\def\ee{\end{equation}}
\def\ea{\end{eqnarray}}

\def\gl{{\mathbf {g}}}
\def\gln{{\mathbf {gl}}_n}
\def\h{{\mathbf{h}}}

\def\LL{\operatorname{L}}
\def\Lprime{{\LL'}}
\def\la{\lambda}

\def\lcd{,\dots,}

\def\M{\operatorname{M}}

\def\n{\mathbf {n}}

\def\P{\operatorname{P}}
\def\PPsi{\hat{\Psi}}
\def\R{\operatorname{R}}
\def\RR{\hat{\R}}

\def\q{\check{\operatorname{q}}}
\def\QQ{\operatorname{Q}}
\newcommand{\s}[2]{s_{#1 #2}}
\newcommand{\sprime}[2]{{s'}^{#1}_{#2}}
\def\scirc{{\scriptstyle{\,\diamond\,}}}
\def\si{\sigma}

\def\T{{\operatorname{T}}}
\def\TT{\hat{\T}}
\def\th{\tilde{h}}
\def\U{\operatorname{U}}
\def\Uh{\bar{\U}(\h)}
\def\Ug{\U(\mathbf{g})}
\def\UUg{\bar{\U}(\mathbf{g})}
\def\UU{\hat{\operatorname{S}}}
\def\ve{\varepsilon}
\def\vphi{\varphi}

\def\Z{\mathbb Z}
\def\Sym{\sf{S}}
\renewcommand{\L}[2]{{{\LL}}^{#1}_{#2}}
\newcommand{\tL}[2]{\tilde{{\LL}}_{#1#2}}
\begin{document}
\begin{center}
{\Large\bf Diagonal reduction algebra and reflection equation}

\vspace{.4cm}
{\bf S. Khoroshkin$^\diamond$$^\star$ and
O. Ogievetsky$^{\circ}$\footnote{On leave of absence from P. N. Lebedev Physical Institute, Leninsky Pr. 53,
117924 Moscow, Russia} }

\vskip .2cm
$\diamond\ ${ITEP,  B.Cheremushkinskaya 25, Moscow 117218, Russia}\\
 $\star\ ${National Research University Higher School of Economics,\\  Myasnitskaya 20, Moscow 101000, Russia}\\
\vspace{.1cm}
$\circ\ ${Aix Marseille Universit\'{e}, Universit\'{e} de Toulon, \\CNRS, CPT UMR 7332, 13288, Marseille, France}
\end{center}

\vskip .2cm
\begin{abstract}
\noindent We describe the diagonal reduction algebra $\DRn$ of the Lie algebra $\gln$ in the $R$-matrix formalism.
As a byproduct we present two families of central elements and the braided 
bialgebra structure of $\DRn$.
\end{abstract}
\section{Introduction}  Reduction algebras were introduced in  \cite{M,AST2} for a study of
representations of a Lie algebra with the help of the restriction to the space of highest weight
 vectors with respect to a reductive subalgebra $\gl$. In an abstract setting, for an associative algebra
$\A$ which contains $\U(\gl)$ as a subalgebra and satisfies certain finiteness conditions,
 the corresponding reduction algebra is the double coset \rf{coset} equipped with a
 nontrivial multiplication \rf{mult}.

This associative multiplication, defined with the help of the extremal projector of Asherova-Smirnov-Tolstoy
\cite{AST}, can be also described \cite{K} by means of the so called universal dynamical twist
 $J$. This twist gives rise to a solution of the universal dynamical
Yang--Baxter equation \cite{ABRR}.

The diagonal reduction algebra $\DR$ is a particular case, associated to the
diagonal embedding of $\Ug$ into $\Ug\otimes\Ug$, of reduction algebras. The algebra $\DR$ acts in the space of highest weight vectors
of the tensor product of two representations of $\gl$, considered as the representation of $\gl$.
 In \cite{KO2,KO3} we  presented a list of ordering defining
relations for natural generators of the diagonal reduction algebra $\DRn$
 of the Lie algebra $\gln$.

The main goal, proposition \ref{prop4}  of this paper, is to relate
the algebra $\DRn$ to the $R$-matrix formalism.
We exhibit a matrix $\LL$ of certain generators
of the algebra $\DRn$ such that the defining relations can be collected
into the operator equation usually called the reflection equation.
  
To this end 
we study and use the  reduction algebras $\Diff(n,N)$ of  algebras of differential operators in $nN$ variables, which we call the algebras of $\h$-deformed differential operators. These reductions were used in \cite{KN} for the representation theory of Yangians. The algebras $\Diff(n,N)$ are closely related, by means of the generalized Harish-Chandra isomorphism \cite{KNV}, to ``relative Yangians" of A. Joseph \cite{J} and ``family algebras" of A. A. Kirillov \cite{Kr}. 
 
 In Section \ref{sechdiff} we introduce a distinguished set of Heisenberg type generators of 
 the algebras  $\Diff(n,N)$ and write down the defining relations for them. It turns out that the relations admit the $\RR$-matrix form, with the well-known solution $\RR$ of the dynamical Yang--Baxter equation. Next we exhibit homomorphisms from the diagonal reduction algebra $\DRn$ to
 $\Diff(n,N)$ ($\h$-analogues of the ``oscillator representations"), which lead to  the presentation of  $\DRn$ by means of the reflection equation.

The description of the algebra $\DRn$ as the reflection equation algebra
has many advantages. As an application we find two natural families
of central elements of $\DRn$ expressed as
 ``quantum'' traces of powers of $L$-operators, see Section \ref{section4.2}.
 Also, the $R$-matrix formalism reveals the braided
bialgebra structure of $\DR$, see Section \ref{section4.3}.

 The modules over the algebra $\DR$ and the braided tensor
 structure on certain categories of $\DR$-modules we study in a forthcoming publication.

 \section{Reduction algebras}

 In this section we recall the definition and some basic properties of reduction algebras. We restrict ourselves
to reduction algebras related to general linear Lie algebra $\gln$ which we denote further by
$\gl$. Let ${{e}}_{ij}$, $i,j=1,\dots, n$, be the standard generators of the Lie algebra $\gl$,
with the commutation relations
\begin{align}[{
{e}}_{ij},{
{e}}_{kl}]=\delta_{jk}{
{e}}_{il}-\delta_{il}{
{e}}_{kj}\ .\end{align} We use the notation $\h$ and $\n_\pm$ for the Cartan and
 two opposite nilpotent subalgebras of  $\gl$;
 $h_{i}$ denotes the element $e_{ii}\in\h$, and  $\ve_i$ the elements  in $\h^*$, so that
$\ve_i(h_j):=\delta_{i,j}$; $h_{ij}:=h_i-h_j\in \h$ and
 $\th_i:=h_i-i$, $\th_{ij}=\th_i-\th_j$ are
the elements of $\U(\h)$.  We define $\Uh$ to be the ring of fractions of the commutative ring
$\U(\h)$ with respect to  the multiplicative set of denominators, generated by the elements
 $(h_{ij}+k)^{-1}$, $k\in\Z$; and $\UUg$ to be  the ring of fractions of  the universal enveloping algebra
$\U(\gl)$ with respect to  the same set of denominators.

\vskip .2cm
\noindent{\bf 1.} Let $\A$ be an associative algebra which contains $\Ug$ as a subalgebra, and the adjoint action of
$\gl$ on $\A$ is locally finite. In particular
 $\A$ is a $\Ug$-bimodule with respect to the multiplication by elements from $\Ug$ on
 the left and on the right. Assume in addition that
 $\A$ is free as the left $\U(\h)$-module and the adjoint action of $\U(\h)$ is semisimple.
 Let $\Ab$ be the localization $\Ab=\A\otimes_{\U(\h)}\Uh$. The  double coset space
\be\label{coset}\Ac:=\n_-\Ab\backslash \Ab/\Ab\n_+\ee
   equipped with a natural associative multiplication  $\scirc$, see e.g. \cite{KO1} for details,
 is usually called the {\it reduction} (double coset) algebra. The multiplication $\scirc$
is described by the prescription
\be\label{mult}x\scirc y=x\P y\,,\ee
where $\P$ is the extremal projector \cite{AST}. The projector $\P$ belongs to a certain extension of
$\UUg$,  satisfies the properties
\begin{align*} x\P&=\P y=0  &&\text{for} \ x\in \n_+,\ y\in\n_-,
\\ \P&=1 \mod \n_-\UUg,\qquad &&\P=1\mod
 \UUg\n_+,\qquad \P^2=\P.
\end{align*}
and can be given by the explicit multiplicative formula \cite{AST}. Alternatively,
 one can find representatives  $\tilde{x}\in\Ab$ and $\tilde{y}\in \Ab$ of coset classes
$x$ and $y$, such that  $\tilde{x}$ belongs to the normalizer of
  the left ideal $\Ab\n_+$ , or $\tilde{y}$ belongs to the normalizer of
 the right ideal $\n_-\Ab$. Then $x\scirc y$ is  the image in the coset space $\Ac$
of the product $\tilde{x}\cdot\tilde{y}$.

  In the sequel the algebra $\A$ is  free as the left $\Ug$-module with respect to the multiplication
by elements of $\Ug$ on the left, and this $\Ug$-module is generated by a linear space $V$,
invariant with respect to the adjoint action (which is supposed to be locally finite) of the Lie algebra $\gl$.
 In this case the reduction algebra
 is the free left $\Uh$-module, generated by $V$, see \cite{Zh,KO1}. Note also that  the construction of the double coset
 space $\Ac$ does not use the multiplication in $\A$ but
 the $\Ug$-bimodule structure on $\A$ only.
We denote sometimes by $:x:$ the image in $\Ac$
 of the element $x\in\Ab$. This notation  is useful
to distinguish between the multiplication $\cdot$ in the algebra $\A$ and the multiplication $\scirc$
 in the reduction algebra $\Ac$.

\vskip .2cm
\noindent{\bf 2.} The Weyl group of the root system of $\gl$ is the symmetric group $\Sym_n$.
 Let
$\si_1\lcd\si_{n-1}$ be the generators of $\Sym_n$; $\si_i$ corresponds to
the permutation $(i,i+1)$.
 The group $\Sym_n$ acts on vector spaces $\h^*$ and $\h$, so that $\si_i(h_j)=
h_{\si_i(j)}$ and $\si_i(\ve_j)=\ve_{\si_i(j)}$.  These actions are related by
 $\la(\si(h))=\si^{-1}(\la)(h)$,
$\si\in\Sym_n$, $h\in\h$ and $\la\in\h^*$.  The action of $\Sym_n$ on $\h$
 extends to the action of a cover $\tilde{\Sym}_n$ of the group $\Sym_n$ by automorphisms of
the Lie algebra $\gl$.
We denote by the same symbols
 $\si_i$ the following  automorphisms of the algebra $\U(\gl)$:
$$\si_i(x):=\mathrm{Ad}_{\exp(e_{i,i+1})}
\mathrm{Ad}_{\exp(-e_{i+1,i})}\mathrm{Ad}_{\exp(e_{i,i+1})}(x)\ ,$$
so that
\begin{align*}\si_i(e_{kl})=(-1)^{\delta_{ik}+\delta_{il}}e_{\sigma_i(k)\sigma_i(l)}\
 .\end{align*}
Let $\rho=-\sum_{k=1}^n k\ve_k$. Then the \textit{shifted action\/}
$\circ$ of the group $\Sym_n$ on the vector space $\h^*$ is defined by
setting
\be\label{shifted}
\si\circ\la:=\si(\la+\rho)-\rho \ .
\ee
With the help of \rf{shifted} we induce the action $\circ$ of $\Sym_n$ on the commutative
algebra $\U(\h)$
by regarding the elements of this algebra as
polynomial functions on $\h^\ast$. We have
\begin{align*}\si\circ\th_k=\th_{\si(k)}\ ,\qquad \si\circ\th_{ij}=\th_{\si(i)\si(j)}\ .\end{align*}

We also use the following notation for shift automorphisms of the rings
$\U(\h)$ and $\Uh$. For any $\alpha\in\h^*$ and any element $f\in\Uh$
 denote by $f[\alpha]$ the
image of $f$ under the shift automorphism  $\Uh\to\Uh$, defined by the rule
 $h\to h+(\alpha,h)$. In particular, $h_i[\ve_j]=h_i+\delta_{ij}$.

\vskip .2cm
\noindent{\bf 3.} Assume now that the action of $\tilde{\Sym}_n$ on $\Ug$ extends
to the action by automorphisms of the algebra $\A$.
 For any $i=1,...,n-1$ define the  map  $\q_i:\A\to \Ac$ by
\be\label{q}\q_i(x)=\sum_{k\geq 0}\frac{(-1)^k}{k!}\hat{e}_{i,i+1}^k\left(\sigma_i(x)
\right)e_{i+1,i}^k\prod_{j=1}^k(\th_{i,i+1}-j)^{-1}
\ee
 Here $ \hat{e}_{i,i+1}(x)=[{e}_{i,i+1},x]$ is the adjoint action of ${e}_{i,i+1}$ on
$x$.  The image of
the sum $\n_-\A+\A\n_+$ of the ideals  is a subspace of  $\n_-\Ab+\Ab\n_+$, see
\cite{Zh}, so the formula \rf{q} defines the map, denoted by the same symbol
 $\q_i:\n_-\A\setminus\A/\A\n_+\to \Ac$. The map $\q_i$
satisfies the relations
\be\label{q2}\q_i(hx)=(\si_i\circ h) \q_i(x),\qquad
\q_i(xh)=\q_i(x)(\si_i\circ h)
\ee
for any $h\in\U(\h)$ and $x\in\A$. The use of \rf{q2} extends the map
 $\q_i:\n_-\A\setminus\A/\A\n_+\to \Ac$
to the map $\q_i:\Ac\to\Ac$. The  maps $\q_i$ satisfy the braid group relations
 \cite{Zh} and are automorphisms of the reduction algebra $\Ac$, see \cite{KO1}.

\section{Algebra of $\h$-deformed differential operators}\label{sechdiff}
\subsection{Reduction of the algebras of differential operators}\label{hdefdiff}
 There is a homomorphism
 $\psi:\Ug\to \Dif (n)$ of the algebra $\Ug$ to the algebra of
polynomial differential operators in $n$ variables $x_1,\ldots, x_n$. The image of  $\E{i}{j}\in \Ug$ is
\begin{equation}\label{O51}
\psi(\E{i}{j})=x^i\partial_j\ .
\end{equation}
Denote by $\Diff(n)$ the reduction algebra  of $\Dif(n)\otimes \Ug$
with respect to the diagonal embedding of $\Ug$. We
call $\Diff(n)$ {\it algebra of $\h$-deformed differential operators}
(the explicit definition, given below, makes it clear that this
algebra admits a well-defined limit $\th_{ij}\to\infty$ for $i<j$
and $\th_1>\th_2>\dots>\th_n$ in which it becomes the usual algebra
of differential operators with polynomial coefficients).

The algebra $\Diff(n)$ is generated over $\Uh$ by the classes of  $x^i$ and $\partial_j$,
which we denote by the same symbols. These elements are subject in $\Diff(n)$ to quadratic-linear
 relations over $\Uh$, which we now describe.

These relations can be computed directly. However, almost all the
required information can be found in the description of the  reduction
algebra of $\U(\mathbf{gl}_{n+1})$ with respect to $\U(\gln)$ (this description is a
 basic step in the derivation of the
 Gelfand--Tsetlin basis in \cite{Zh}).  Indeed,
the generators $\E{i,}{n+1}$ and $\E{n+1,}{i}$, $i=1,\ldots, n$, of the corresponding
 double coset algebra form the bases of the fundamental representation
 $\omega$ and its dual $\omega^*$ with
respect to the adjoint representation of $\gl=\gln$. The elements
$x^i$ and $\partial_i$ form the same bases of $\omega$ and
$\omega^*$, the only difference is that the commutators
$[\E{i,}{n+1},\E{n+1,}{j}]$ belong to the Cartan subalgebra of
$\mathbf{gl}_{n+1}$ while $[x^i,\partial_j]=\delta_{ij}$.
  Thus, by \cite[4.5.3]{Zh}, we have
\be \label{O5}\begin{split}
x^i\scirc x^j&=\alpha_{ij}x^j\scirc x^i, \qquad \partial_j\scirc\partial_i=
\alpha_{ij}\partial_i\scirc\partial_j, \qquad i<j,\\
x^i\scirc\partial_j&=\partial_j\scirc x^i\qquad i\not=j,\\
x^i\scirc\partial_i&=\sum_j\beta_{ij}\partial_j\scirc x^j+\mu_i.
\end{split}
\ee
 Here
\be\label{O6} \alpha_{ij}=\frac{\th_{ij}+1}{\th_{ij}},\qquad
\beta_{ij}=\frac{1}{1-\th_{ij}}\frac{\varphi_{j}[\ve_j]}{\varphi_i}\quad\text{with}
\quad
 \varphi_j=\prod_{k: k>j} \frac{\th_{jk}}{\th_{jk}-1},
\ee
and $\mu_i$ are the elements of $\Uh$ which are to be determined by
 another argument.
\begin{lem}\label{lem1}
 We have
$$\mu_i=-\varphi_i^{-1}\ .$$
\end{lem}

\noindent {\it Proof}. Note first that $\mu_n=-1$. Indeed, since $\partial_n$ is a
highest weight vector with respect to the adjoint action of $\Ug$ on $\Dif(n)$,
 $$x^n\scirc\partial_n=:x^n\partial_n:=:\partial_nx^n:-1.$$ On the other hand,  the products
$\partial_j\scirc x^j$ are equal to  sums $\sum_m a_m :\partial_m x^m:$ with some $a_m\in\Uh$ and do not contain
a constant term. Thus $\mu_n=-1=-\varphi_n^{-1}$. For the derivation of other $\mu_j$ we
use the Zhelobenko automorphisms  $\q_i$.

 It is not difficult to check that
\be\label{O7}\begin{split}
\q_i(x^i)&=-x^{i+1}\frac{\th_{i,i+1}}{\th_{i,i+1}-1},
\qquad \q_i(x^{i+1})=x^{i},\qquad \q_i(x^j)=x^j,\quad j\not=i,i+1,\\[-.2em]
\q_i(\partial_{i+1})&=\partial_i\frac{\th_{i,i+1}}{\th_{i,i+1}-1},\qquad
\q_i(\partial_{i})=-\partial_{i+1},\quad \q_i(\partial_j)=\partial_j,\quad
j\not=i,i+1.
\end{split}
\ee
 Now we apply the automorphism $\q_{n-1}$ to the  already known identity
$$x^n\scirc\partial_n=\sum_j\beta_{n,j}\partial_j\scirc x^j-1.$$
We get the relation
$$\frac{\th_{n-1,n}}{\th_{n-1,n}-1}x^{n-1}\scirc\partial_{n-1}
=\sum_j\q_{n-1}\left(\beta_{n,j}\partial_j\scirc x^j\right)-1.$$
Taking into account the last line in (\ref{O5}), we rewrite this relation in the form
$$x^{n-1}\scirc\partial_{n-1}=\sum_j\beta_{n-1,j}\partial_j\scirc x^j-
\frac{\th_{n-1,n}-1}{\th_{n-1,n}},$$
 which  implies that
$$\mu_{n-1}=-\frac{\th_{n-1,n}-1}{\th_{n-1,n}}=-\vphi_{n-1}^{-1}.$$
 Next we apply $\q_{n-2}$ to the identity
$x^{n-1}\scirc\partial_{n-1}=\sum_j\beta_{n-1,j}\partial_j\scirc x^j-\vphi_{n-1}^{-1}$
 and find that $\mu_{n-2}=-\vphi_{n-1}^{-1}$
 and then further
$\mu_i=-\vphi_i^{-1}$ for all $i$. \hfill $\square$

\vskip .4cm
 The last line in \rf{O5} can be now rewritten as
$$x^i\scirc\partial_i=\sum_j\frac{1}{1-\th_{ij}}\,\partial_j\vphi_j\scirc x^j
\vphi_i^{-1}-\vphi_i^{-1},$$
which suggests the  change of variables
\be\label{O8} \der_j:=\partial_j \vphi_j.
\ee
In the new variables the relations \rf{O5} and \rf{O7} look as follows:
\begin{equation}\label{O9}
\begin{split}
x^i\scirc x^j&=\frac{\th_{ij}+1}{\th_{ij}}x^j\scirc x^i,\ i<j,\quad
\der_i\scirc\der_j=\frac{\th_{ij}-1}{\th_{ij}}\,\der_j\scirc \der_i,\qquad i<j,\\
x^i\scirc \der_j&=\der_j\scirc x^i, \ i<j;\qquad \qquad
x^i\scirc \der_j=\frac{\th_{ij}(\th_{ij}-2)}{(\th_{ij}-1)^2}\,\der_j\scirc x^i,
\ i>j,\\
x^i\scirc \der_i&= \sum_j\frac{1}{1-\th_{ij}}\,\der_j\scirc x^j-1,
\end{split}
\end{equation}
\be\label{O7a}\begin{split}
\q_i(x^i)&=-x^{i+1}\frac{\th_{i,i+1}}{\th_{i,i+1}-1},
\qquad \q_i(x^{i+1})=x^{i},\quad \q_i(x^j)=x^j,\quad j\not=i,i+1,\\
\q_i(\der_{i})&=- \frac{\th_{i,i+1}-1}{\th_{i,i+1}}\der_{i+1},\qquad
\q_i(\der_{i+1})=\der_{i},\quad \q_i(\der_j)=\der_j,\quad
j\not=i,i+1.
\end{split}
\ee
We have other sets of Heisenberg type generators in the algebra $\Diff(n)$.
Set
\begin{equation}
  \vphi'_j=\prod_{k: k<j} \frac{\th_{jk}}{\th_{jk}-1},
\end{equation}
and put
\be\label{O8a}\bar{\der}_j:= \partial_j {\vphi'_j}^{-1}.
\ee
Arguments parallel to that of Lemma \ref{lem1} show that
\begin{equation}\label{O9a}
\begin{split}
x^i\scirc x^j&=\frac{\th_{ij}+1}{\th_{ij}}x^j\scirc x^i,  i<j,\qquad
\bar{\der}_i\scirc\bar\der_j=\frac{\th_{ij}-1}{\th_{ij}}\,\bar\der_j\scirc \bar\der_i,\qquad i<j,\\
 \bar\der_j\scirc x^i&=x^i\scirc\bar\der_j, \ i>j,\qquad
 \bar\der_j\scirc x^i=\frac{\th_{ij}(\th_{ij}-2)}{(\th_{ij}-1)^2}\,x^i\scirc\bar\der_j,
\ i<j,\\
\bar\der_i\scirc x^i &= \sum_j\frac{1}{1+\th_{ij}}\,x^j\scirc\bar\der_j+1.
\end{split}
\end{equation}
and the same as in \rf{O7a} coefficients in the action of Zhelobenko
automorphisms. Alternatively, we can leave the variables
$\partial_j$ unchanged and rescale the variables $x^j$. For the
variables $\bar{x}^j:=\vphi_j x^j$ and $\partial_j$ we get the
relations of the form \rf{O9} with the reversed inequalities between
$i$ and $j$ in first two lines of the relations \rf{O9}; for the
variables $\bar{\bar{x}}^j:={\vphi'_j}^{-1}x^j$ and $\partial_j$ we
get the relations of the form \rf{O9a} with the reversed
inequalities between $i$ and $j$ in first two lines of relations
\rf{O9a}.

\subsection{Polarized form  of relations}\label{polfore}
\paragraph{1.} The first line of \rf{O9} can be written in the following 
polarized form:
 \be
\begin{split}\label{O10}
x^{i}\scirc x^{j}&=\frac{1}{\th_{ij}}x^{i}\scirc x^{j}+
\frac{\th_{ij}^2-1}{\th_{ij}^2}x^{j}\scirc x^{i},\qquad i<j,\\[-.4em]
x^{i}\scirc x^{j}&=\frac{1}{\th_{ij}}x^{i}\scirc x^{j}+
x^{j}\scirc x^{i},\qquad \qquad i>j\, ,
\end{split}
\ee
\be
\begin{split}\label{O10a}
\der_{i}\scirc \der_{j}&=-\frac{1}{\th_{ij}}\der_{i}\scirc \der_{j}+
\frac{\th_{ij}^2-1}{\th_{ij}^2}\der_{j}\scirc \der_{i},\qquad i<j,\\
\der_{i}\scirc \der_{j}&=-\frac{1}{\th_{ij}}\der_{i}\scirc \der_{j}+
\der_{j}\scirc \der_{i},\qquad \qquad i>j\, .
\end{split}
\ee Rewrite \rf{O10}, \rf{O10a} and the last two lines in \rf{O9}
 in an operator form\footnote{Unless the opposite is stated, we adopt the Einstein convention:
if a tensor index  appears in an expression twice, once as an upper
index and once as a lower index, the summation over this index is
assumed.}
\be\label{O11} x^i\scirc x^j=\RR_{kl}^{ij}x^k\scirc
x^l,\qquad \der_i\scirc \der_j=\RR_{ji}^{lk}\der_k\scirc
\der_l,\qquad
 x^i\scirc\der_j=\TT_{jl}^{ik}\der_k\scirc x^l-\delta_{j}^i,
\ee where  $\RR_{kl}^{ij}$
 and $\TT_{jl}^{ik}$ are matrix
coefficients of operators
 $\RR,\TT: \C^n\otimes\C^n\to\Uh\otimes\C^n\otimes\C^n$. Their nonzero values are
\be\label{O12}
\RR_{ij}^{ij}=\frac{1}{\th_{ij}},\ i\not=j,\qquad \RR_{ji}^{ij}=\left\{
\begin{array}{cc}
 \dfrac{\th_{ij}^2-1}{\th_{ij}^2},& \, i<j,\\[0.5em]
 1,&\, i\geq j
\end{array}
\right. \ee \be\label{O14} \TT_{ij}^{ij}=-\frac{1}{\th_{ij}-1},\
i\not=j,\qquad \TT_{ij}^{ji}=\left\{
\begin{array}{cc}
 \dfrac{\th_{ij}(\th_{ij}+2)}{(\th_{ij}+1)^2},& \, i<j,\\[.5em]
 1,&\, i\geq j
\end{array}
\right. \ee Note the following identity:
\be\label{O15}
\TT_{ij}^{kl}[-\ve_l]=\RR_{ji}^{lk}. \ee
Similarly, the relations
\rf{O9a} can be presented in an operator form as \be\label{Q3}
x^i\scirc x^j=\RR_{kl}^{ij}x^k\scirc x^l,\qquad \bar\der_i\scirc
\bar\der_j=\RR_{ji}^{lk}\bar\der_k\scirc \bar\der_l,\qquad
 \bar\der_j\scirc x^i=\UU_{jl}^{ik}x^l\scirc \bar\der_k+\delta_{j}^i,
\ee
where the non-zero elements of $\UU$ are
\be\label{Q4}\UU^{ij}_{ij}=\frac{1}{\th_{ij}+1}\ ,\qquad \UU^{ij}_{ji}=\left\{\begin{array}{cc}
1&,\ i>j\\[.5em]
\displaystyle{\frac{\th_{ij}(\th_{ij}-2)}{(\th_{ij}-1)^2} }&,\ i<j\end{array}\right.
\ee

\paragraph{2.} The relations in the last two lines of \rf{O9a}
can be obtained in another way, by inverting the last relation in \rf{O11}.
 Let $\PPsi$ be the skew inverse to $T$ (see e.g. \cite{O}, section 4.1.2 for details of the $R$-matrix technique needed here), that is,
\be\label{skew}
\PPsi_{jl}^{ik}\TT_{kn}^{lm}=\delta_{n}^i\delta_j^m.
\ee

Multiplying the relation
$$x^l\scirc\der_k=\TT_{kn}^{lm}\der_m\scirc x^n-\delta_k^l$$
by $\PPsi_{jl}^{ik}$ from the left 
and contracting repeated indices,
 we get the relation
\be\label{Qdop1}\der_j \scirc x^i=\PPsi_{jl}^{ik}x^l\scirc \der_k +\PPsi_{jk}^{ik}.\ee
The relations \rf{O8}, \rf{O8a}, and the last two lines of \rf{O9a} imply that
\be\label{Q1}\bar{\partial}_i=\bar{\bar{\partial}}_i (\QQ^-_i)^{-1} = \QQ^+_i\bar{\bar{\partial}}_i,\ee
where
\be \label{Qpm} \QQ^{\pm}_i=\prod_{k:k\neq i}\frac{\th_{ik}\pm 1}{\th_{ik}},\qquad \QQ^-_j[\ve_j]\QQ^+_j=1\,,\ee
and
\be\label{explpsi}\PPsi^{ij}_{ij}=\QQ^+_i\QQ^-_j\frac{1}{\th_{ij}+1}\ ,\qquad \PPsi^{ij}_{ji}=\left\{\begin{array}{cc}
1&,\ i<j\\[.5em]
\displaystyle{\frac{(\th_{ij}-1)^2}{\th_{ij}(\th_{ij}-2)} }&,\ i>j\end{array}\right.\ee
Comparing \rf{Qdop1} and the last line in \rf{O9a} we conclude that
\be\label{Q2}\PPsi_{jk}^{ik}=\QQ^+_i\delta_j^i .\ee
Define operators $\QQ^\pm:\C^n\to\C^n$ by
$$(\QQ^\pm)_j^i=\QQ^\pm_j\delta_j^i.$$
Then the relation \rf{Q2} can be rewritten as
$$\text{Tr}_2 \PPsi_{12}=\QQ^+_1.$$
Here the lower index specifies the number of the copy of the space $\C^n$ in which the corresponding
operator acts nontrivially. For example, $\QQ^+_1$ stands for the operator $\QQ^+$ acting in the first copy and
$\text{Tr}_2$ means the contraction of indices in the second copy.

\subsection{Reflection equation and copies}

To lighten the notation, in the formulation of statements, the matrix multiplication of matrices with entries in a reduction algebra is written without the symbol $\scirc$
(which is always assumed). 

Set $\tilde{\text{L}}^{i}_{j}:=x^i\scirc\der_j\in\Diff(n)$. 
\begin{prop} \label{prop3.2} The matrix $\tL{}{}$ satisfies the reflection equation
\be\label{O15a}\RR_{12}\tL_1\RR_{12}\tL_1-\tL_1\RR_{12}\tL_1\RR_{12}=
\RR_{12}\tL_1-\tL_1\RR_{12}.\ee 
 \end{prop}
\renewcommand{\s}[2]{s^{#1}_{ #2}}
\renewcommand{\L}[2]{{{\LL}}^{#1}_{#2}}
\renewcommand{\tL}[2]{\tilde{{\LL}}^{#1}_{#2}}
\noindent{\it Proof}.  Consider the monomial $x^{i_1}\scirc x^{{i_2}}\scirc \der_{{j_1}}\scirc\der_{{j_2}}$.
 Reorder it in two ways.  The first way:
\begin{align*}
 & x^{{i_1}}\scirc x^{{i_2}}\scirc \der_{{j_1}}\scirc\der_{{j_2}}=\RR^{{i_1}{{i_2}}}_{kl}x^k\scirc x^l
\scirc \der_{{j_1}}\scirc\der_{{j_2}}=\RR^{{i_1}{{i_2}}}_{kl}x^k\scirc\TT^{l{n}}_{{{j_1}} m}\der_{n}
\scirc x^m\scirc\der_{{j_2}}-\RR^{{i_1}{{i_2}}}_{kl}\delta_{{j_1}}^l x^k\scirc \der_{{j_2}}\\[.4em]  =
 &\sum_{k,l,m,n} \RR^{{i_1}{{i_2}}}_{kl}x^k\circ\der_{n}\scirc\TT^{l{n}}_{{{j_1}} m}[-\ve_{n}]
x^m\scirc\der_{{j_2}} -\RR^{{i_1}{{i_2}}}_{k{{j_1}}} x^k\scirc \der_{{j_2}}=
\RR^{{i_1}{{i_2}}}_{kl} \tL{k}{{n}}\scirc\RR^{{n} l}_{{{j_1}m}}\tL{m}{{{j_2}}}-
\RR^{{i_1}{{i_2}}}_{k{{j_1}}}\tL{k}{{{j_2}}}\,.
\end{align*}
The second way  of reordering:
\begin{align*}
&x^{{i_1}}\scirc x^{{i_2}}\scirc \der_{{j_1}}\scirc\der_{{j_2}}=x^{{i_1}}\scirc x^{{i_2}}\scirc
\RR^{{m}{n}}_{{{j_2}}{{j_1}}}\der_{n}\scirc\der_{m}\\[.4em]
=&
\RR^{{m}{n}}_{{{j_2}}{{j_1}}}[-\ve_{i_1}-\ve_{{i_2}}]x^{{i_1}}\scirc\TT^{{{i_2}}{t}}_{{n} k}
\der_{t}\scirc x^k\scirc\der_{m}-
\RR^{{m}{n}}_{{{j_2}}{{j_1}}}[-\ve_{i_1}-\ve_{{i_2}}]{\delta}^{{i_2}}_{n} x^{{i_1}}\scirc \der_{m} \\[.4em]  =
&\RR^{{m}{n}}_{{{j_2}}{{j_1}}}[-\ve_{i_1}-\ve_{{i_2}}]
\TT^{{{i_2}}{t}}_{{n} k}[-\ve_{i_1}]x^{{i_1}}\scirc\der_{t}\scirc x^k\scirc\der_{m}
- \RR^{{m}{i_2}}_{{{j_2}}{{j_1}}}[-\ve_{i_1}-\ve_{i_2}] x^{{i_1}}\scirc \der_{m} \\[.4em]  =
&\sum_{t,k,m,n}x^{{i_1}}\scirc\der_{t}\scirc\TT^{{{i_2}}{t}}_{{n} k}[-\ve_{t}]
x^k\scirc\der_{m}\RR^{{m}{n}}_{{{j_2}}{{j_1}}}
[-\ve_{t}-\ve_{{i_2}}+\ve_k-\ve_{m}] -
\sum_{m} x^{{i_1}}\scirc \der_{m}
\RR^{{m}{i_2}}_{{{j_2}}{{j_1}}}[-\ve_{m}-\ve_{i_2}] \\[.4em] = &
\sum_{t,k,m,n}\tL{i}{{t}}\TT^{{{i_2}}{t}}_{{n} k}[-\ve_{t}]\scirc\tL{k}{{m}}\RR^{{m}{n}}_{{{j_2}}{{j_1}}}
[-\ve_{t}-\ve_{{i_2}}+\ve_k-\ve_{m}]-
\sum_{m} \tL{i_1}{{m}}\RR^{{m}{i_2}}_{{{j_2}}{{j_1}}}[-\ve_{m}-\ve_{i_2}]\,.
\end{align*}
 Matrix elements  $\RR^{{i}{j}}_{{{k}}{{l}}}$  and $\TT^{{i}{j}}_{{{k}}{{l}}}$ are nonzero only if $i=k$ and $j=l$, or $i=l$ and $j=k$. We can thus replace the shift by
$\ve_{t}+\ve_{{i_2}}$ in the last displayed line with the shift by $\ve_{n}+\ve_k$.
Besides,  the matrix coefficient $ \RR^{{m}{n}}_{{{j_2}}{{j_1}}}$ depends only
 on the difference $h_{n}-h_{m}$ and is invariant with respect to the shift
by $\ve_{m}+\ve_{n}$. Using the relation \rf{O15} we rewrite the
 result as
\be\label{O17}
x^{{i_1}}\scirc x^{{i_2}}\scirc \der_{{j_1}}\scirc\der_{{j_2}}=
\tL{i_1}{{t}}\RR^{ {{i_2}}{t}}_{k {n} }\scirc\tL{k}{{m}}\RR^{{m}{n}}_{{{j_2}}{{j_1}}}
- \tL{i_1}{{m}}\RR^{{m}{i_2}}_{{{j_2}}{{j_1}}}
\ee
Comparing these two ways, we obtain the equality
$$
\RR^{ij}_{kl} \tL{k}{{n}}\scirc\RR^{{n} l}_{m{s}}\tL{m}{{t}}-
\RR^{ij}_{k{s}}\tL{k}{{t}}=
\tL{i}{{t}}\RR^{{t} j}_{k {n} }\scirc\tL{k}{{u}}\RR^{{u}{n}}_{{t}{s}}
- \tL{i}{{u}}\RR^{{u}{n}}_{{t}{s}},$$
that is the relation \rf{O15a}. \hfill $\square$
\bigskip

Let now $\Dif(n,N)$ be the algebra of differential operators in $nN$ variables $x^{i\alpha}$,
$i=1,\ldots,n$, $\alpha=1,\ldots, N$.
 There is a homomorphism
 $\psi_N:\Ug\to \Dif(n,N)$.  The image of  $\E{i}{j}\in \Ug$ is
\begin{equation}\label{O32}\psi(\E{i}{j})=\sum_\alpha x^{i\alpha}\partial_{j\alpha}.
\end{equation}

Denote by $\Diff(n,N)$ the reduction algebra  of $\Dif(n,N)\otimes \Ug$ with respect to the diagonal embedding
 of $\Ug$. The algebra $\Diff(n,N)$ is generated over $\Uh$ by the classes of all
$x^{i\alpha}$ and $\partial_{j\beta}$, which we denote by the same symbols.
  Denote
\be\label{O36}
\der_{j\beta}=\partial_{j\beta}\varphi_j\,.\ee
The calculations from Sections \ref{hdefdiff} and \ref{polfore} can be repeated for any $N$. The result is
\begin{prop} \label{prop2} The elements $x^{i\alpha}$ and $\der_{j\beta}$
satisfy the following relations
 \begin{equation}
\begin{split}
  \label{O31}
x^{i\alpha}\scirc x^{j\beta}&=\RR_{kl}^{ij}x^{k\beta}\scirc x^{l\alpha},\qquad
\der_{i\alpha}\scirc \der_{j\beta}=\RR_{ji}^{lk}\der_{k\beta}\scirc \der_{l\alpha},\\
 &x^{i\alpha}\scirc\der_{j\beta}=\RR_{lj}^{ki}[\ve_l]\der_{k\beta}\scirc x^{l\alpha}-
\delta^\alpha_\beta\delta_{j}^{i}\, .
\end{split}
 \end{equation}
\end{prop}

\newcommand{\teL}{\operatorname{L}}

The same proof as that of proposition \ref{prop3.2} shows that the combinations
\begin{equation}\label{O33}\tL{i}{j}:=\sum_\alpha x^{i\alpha}\scirc\der_{j\alpha}
\end{equation}
 satisfy the same reflection equation. We formulate this assertion in the separate proposition.
\begin{prop} \label{prop3} The matrix $\tL{}{}$, defined in \rf{O33} satisfies the reflection equation
\be\label{O35a}\RR_{12}\tilde{\teL}_1\RR_{12} \tilde{\teL}_1-
\tilde{\teL}_1\RR_{12}\tilde{\teL}_1\RR_{12}=\RR_{12}\tilde{\teL}_1-
\tilde{\teL}_1\RR_{12}.\ee
 \end{prop}
{\it Remark.} Equally well one can define the grassmanian version of $\h$-deformed differential
operators, the reduction of the algebra $\Ug\otimes\C[\xi^{i\alpha},d_{j\beta}]$ where
$ \xi^{i\alpha}$ are anti-commuting variables and $d_{j\beta}$ grassmanian derivatives in them.
The defining relations differ by the following sign changes:
\begin{equation}
\begin{split}
\xi^{i\alpha}\scirc \xi^{j\beta}&=-\RR_{kl}^{ij}\xi^{k\beta}\scirc \xi^{l\alpha},\qquad
\bar d_{i\alpha}\scirc \bar d_{j\beta}=-\RR_{ji}^{lk}\bar d_{k\beta}\scirc \bar d_{l\alpha},\\
 &\xi^{i\alpha}\scirc\bar d_{j\beta}=-\RR_{lj}^{ki}[\ve_l]\bar d_{k\beta}\scirc xi^{l\alpha}-
\delta^\alpha_\beta\delta_{j}^i\, .
\end{split}
 \end{equation}
Here $\bar d_{j\beta}:=d_{j\beta}\varphi_j$. Then the operator $\tL{}{}$ with entries
$$\tL{i}{j}:=\sum_\alpha \xi^{i\alpha}\scirc\bar d_{j\alpha}$$
 satisfies the same reflection equation \rf{O35a}.

\subsection{$R$-matrix and its skew inverse}
Using the first relation in (\ref{O31}) we  can reorder a monomial $x^{i\alpha}\scirc x^{j\beta}\scirc x^{k\gamma}$, $\alpha\not=\beta\not=\gamma \not=\alpha$, as a 
combination of monomials of the form $x^{\bullet\gamma}\scirc x^{\bullet\beta}\scirc x^{\bullet\alpha}$
in two ways, as $(x^{i\alpha}\scirc x^{j\beta})\scirc x^{k\gamma}$ or as $x^{i\alpha}\scirc (x^{j\beta}\scirc x^{k\gamma})$. The ordered products form a basis in the reduction algebra 
(see \cite{KO2}, section 2) so the two ways of reordering lead to the same result. This implies certain compatibility conditions for the operator $\RR$ which are 
formulated in the following Proposition. 

\begin{prop}
The operator $\RR$ is a solution of the dynamical Yang--Baxter equation
\be \sum_{a,b,u}{\RR}^{ij}_{ab}\RR^{bk}_{ur}[-\varepsilon_a]\RR^{au}_{mn}=\sum_{a,b,u}
\RR^{jk}_{ab}[-\varepsilon_i]\RR^{ia}_{mu}
\RR^{ub}_{nr}[-\varepsilon_m]\ .\ee
\end{prop}
This solution has already appeared several times in different contexts (see e.g. \cite{I,ES} and references therein).

In addition, the operator $\RR$ satisfies the relations
\be\label{Q5do}\RR^2=\text{Id}_{\C^n\otimes\C^n}\, ,\ee
\be\label{Q5}\RR_{21}=\RR^T\vert _{\th\mapsto -\th},\ee
where $(\RR_{21})_{jl}^{ik}:=\RR_{lj}^{ki}$ and $(\RR^{T})_{jl}^{ik}:=\RR_{ik}^{jl}$ and
$$\QQ^+_i[-\epsilon_i]\,\QQ^+_k[-\epsilon_i-\epsilon_k]\,
\RR^{ik}_{jl}=\RR^{ik}_{jl}\,\QQ^+_j[-\epsilon_j]\,\QQ^+_l[-\epsilon_j-\epsilon_l]\ ,$$
which is an immediate consequence of
$$\QQ^-_j \,\QQ^-_i[-\epsilon_j]=\QQ^-_i \,\QQ^-_j[-\epsilon_i]\ .$$
The operators $\PPsi$, $\UU $ and $\RR$ are related by
\be\label{skdop1}\PPsi^{ik}_{jl}=\QQ^+_i[\epsilon_k-\epsilon_l]\,\UU^{ik}_{jl}\,(\QQ^+_l)^{-1}[-\epsilon_l]\ ,\ee
\be\label{skdop2}\UU^{ij}_{kl}=\RR^{ij}_{kl}[\epsilon_k]\ .\ee
By \rf{Q5}, $\PPsi_{21}=\PPsi^T\vert _{\th\mapsto -\th}$ so
\be\label{Q6}\PPsi^{am}_{an}=(\PPsi_{21})^{ma}_{na}=(\PPsi^T\vert _{\th\mapsto -\th})^{ma}_{na}=
(\PPsi\vert _{\th\mapsto -\th})_{ma}^{na}=\QQ^+_n\vert _{\th\mapsto -\th}\delta^n_m=\QQ^-_n
\delta^m_n\ ,\ee
or
$$\text{Tr}_1\PPsi_{12}=\QQ^-_2.$$
It follows from \rf{skew} together with \rf{Q2} and \rf{Q6} that
\be\label{zadsled}\sum_{a} \QQ^-_a[-\epsilon_m]\RR^{ma}_{na}=\delta^m_n\ ,\ee
$$\sum_{a} \QQ^+_a[\epsilon_m]\RR^{am}_{an}=\delta^m_n\ .$$
{\it Remark.} The reordering of the monomial 
$$(\der_{i\alpha}\scirc \der_{j\beta})\scirc \der_{k\gamma}=\der_{i\alpha}\scirc (\der_{j\beta}\scirc \der_{k\gamma})$$
leads to a compatibility condition for the operator $\RR$ which is equivalent to the same 
dynamical Yang--Baxter equation for $\RR$. 

The reordering of any of the monomials 
 \begin{align*}
(x^{i\alpha}\scirc x^{j\beta})\scirc \der_{k\gamma}&=x^{i\alpha}\scirc (x^{j\beta}\scirc \der_{k\gamma}),\\
(x^{i\alpha}\scirc \der_{j\beta})\scirc \der_{k\gamma}&=x^{i\alpha}\scirc (\der_{j\beta}\scirc \der_{k\gamma}),
\end{align*}  
$\alpha\not=\beta\not=\gamma \not=\alpha$, 
similarly leads to a compatibility condition for the operator $\RR$ which (in each case) now is equivalent to the 
dynamical Yang--Baxter equation for $\RR$ together with the equality \rf{Q5do}. The verification of this statement uses 
the fact that the matrix element $\R^{ij}_{kl}$ can be nonzero only if   $i=k$ and $j=l$, or $i=l$ and $j=k$.

\section{Diagonal reduction algebra}
\subsection{$R$-matrix presentation}
The diagonal reduction algebra $\DR$ is by definition the reduction algebra of $\A=\U(\gl)\otimes\U(\gl)$
 with respect to the diagonal embedding of  $\U(\gl)$.
Let $\e{i}{j}^{(1)}$ and $\e{i}{j}^{(2)}$ be the standard generators $e_{ij}$
 of the Lie algebra $\gl$ in the first and the second tensor components.
Denote by $\s{i}{j}$  the generators of the diagonal
reduction algebra defined as the images in $\DR$ of $\e{i}{j}^{(1)}$. In other words,
$$\s{i}{j}=\P\ \e{i}{j}^{(1)}\P.$$
 We will also need another set of generators
$$\sprime{i}{j}=\P\ \e{i}{j}^{(2)}\P.$$
The elements $\s{i}{j}$ and $\sprime{i}{j}$ are related by
\be\label{relsij}
\s{i}{j}+\sprime{i}{j}=h_i\delta^i_j\,.\ee
In addition to the elements $\varphi_j$, defined in \rf{O6}, we need,
for any $j=1,\ldots, n$, and $m>j$,    the following elements of $\Uh$:
$$
\varphi_{jm}=\prod_{k: j<k<m} \frac{\th_{jk}}{\th_{jk}-1}. $$
The description of the algebra $\DR$ in terms of generators $\s{i}{j}-\sprime{i}{j}$
was given in \cite{KO2,KO3}.
Here we suggest another presentation. Let $\L{i}{j}$, $i,j=1,\ldots,n$, be the following elements of $\DR$:
\be\label{Os3}
\L{i}{j}:=\left\{
\begin{array}{cc} \s{i}{j}\varphi_{j},& i\not=j\\
    \left(\s{i}{j}-\sum\limits_{m:m>i}\s{m}{m}\dfrac{1}{\th_{im}
        \varphi_{im}}\right)\varphi_j,
    & i=j
\end{array}
\right.
\ee
The elements $\L{i}{j}$ are linear combinations of $\s{k}{l}$ with the triangular transition matrix.
 Thus $\L{i}{j}$ generate $\DR$ as the  algebra over $\Uh$.

 \begin{prop}\label{prop4}
 Elements $\L{i}{j}\in \DR$ satisfy quadratic-linear relations collected in reflection equation
 \be\label{O35b}\RR_{12}{\teL}_1\RR_{12} {\teL}_1-
 {\teL}_1\RR_{12} {\teL}_1\RR_{12}=\RR_{12}{\teL}_1-
 {\teL}_1\RR_{12}.\ee
The relations \rf{O35b} form a complete list of defining relations
over the field of fractions of $\Uh$.
 \end{prop}
It is plausible that  the relations \rf{O35b} form a complete list of defining relations over the ring $\Uh$ itself.

\vskip .1cm
Proof of proposition \ref{prop4} is based on the properties of  homomorphisms of the diagonal reduction algebra $\DR$
to $ \Diff(n,N)$ and the  proposition \ref{prop3}.

  Using the map  $\psi$, see \rf{O32}, we define a homomorphism
  \begin{equation}
  \label{O52}
\psi_1:=\psi\otimes 1:\ \Ug\otimes\Ug\to\Dif(n,N)\otimes\Ug.
  \end{equation}
 The map $\psi_1$ sends the diagonal $\gl$ to the diagonal $\gl$ and thus defines the homomorphism
 $$\tilde{\psi}_1:\ \DR\to
 \Diff(n,N).$$
 \begin{lem}\label{lem2} The map $\tilde{\psi}_1$ sends the generator $\L{i}{j}\in\DR$ to
    the element $\tilde{\LL}^{i}_{j}\in
    \Diff(n,N)$, (see\rf{O33})
    $$\tilde{\psi}_1(\L{i}{j})=\tilde{\LL}^{i}_{j}.$$
 \end{lem}
 {\it Proof} of proposition \ref{prop4}. The statement of proposition \ref{prop4} follows from proposition \ref{prop3},
lemma \ref{lem2} and the injectivity
 of the map ${\psi}$ for $N\geq n$ (consider the map ${\psi}$ as the tangent map for the group action 
 $\operatorname{GL}_n\times \operatorname{Mat}_{n\times N}\to  \operatorname{Mat}_{n\times N}$;
 for $N=n$ we obtain the classical isomorphism of $\Ug$ with the ring of right invariant differential operators on $\operatorname{GL}_n\subset\operatorname{Mat}_{n\times n}$).  

 In particular, for $N\geq n$ the map $\tilde{\psi}_1$ is injective on the subspace of $\DR$, generated by polynomials in $\L{i}{j}$ of degree not bigger
than two. This proves that the relations \rf{O35b} are satisfied. In  \cite{KO2}, we  presented a complete list
of  relations for the generators $\s{i}{j}-\sprime{i}{j}$  of the algebra $\DR$. These  are ordering
relations of degree  at most two in these generators.
This means, in particular, that the number of linearly independent over $\Uh$ quadratic monomials in these generators equals the number
 of ordered quadratic monomials in $n^2$ generators. Since the transition to the generators $\L{i}{j}$
 has a triangular form,  the relations are of degree at most two on generators $\L{i}{j}$ and the number of linearly
independent quadratic monomials in $\L{i}{j}$ is as before.
 Therefore it is left to prove the completeness of relations \rf{O35b} in degree at most two. To show this we
write the relations  \rf{O35b} in the form
$$Z^Ac^\mu_A(\th)=\text{linear in $\L{i}{j}$ terms},$$
 where $\{Z^A\}$ is the set of quadratic monomials in $\L{i}{j}$.

In the asymptotic regime $\th_1>\th_2>\dots>\th_n$ and $\th_{ij}\to\infty$ for $i<j$ these relations are well defined
and become
 $$Z^Ac^\mu_A(\infty)=\text{linear in $\L{i}{j}$ terms}.$$
Indeed, in this limit the matrix $\RR_{12}$ turns into the permutation matrix $\P_{12}$ so the relations \rf{O35b} turn into
the standard defining relations for the Lie algebra $\gl$,
$$\Lambda_2\Lambda_1-\Lambda_1\Lambda_2 =\P_{12}\left(\Lambda_1-\Lambda_2\right).$$
By general deformation arguments, $\text{rk}\, c^\mu_A(\th)\geq \text{rk}\, c^\mu_A(\infty)$. Hence the number
 of linearly independent quadratic monomials for generic $\th$ is not bigger than that in the limit. But in both cases,
 for generic $\th$ and asymptotically, this number equals the number of ordered (for any linear order)
quadratic monomials in $\L{i}{j}$.
 \hfill{$\square$}

 \medskip
  \noindent{\it Proof} of lemma \ref{lem2}. Denote by $:x^{i\alpha}\partial_{j\alpha}:$  the image of
the element  $x^{i\alpha}\partial_{j\alpha}\otimes 1\in\Dif(n,N)\otimes\Ug$ in the reduction algebra 
$\Diff(n,N)$. Due to the definitions
 of elements $\L{i}{j}$, $\tilde{\LL}^{i}_{j}$ and of the map $\psi_1$ it is sufficient to establish the equalities
 \be\label{O34}
x^{i\alpha}\scirc \der_{j\beta}=\left\{
\begin{array}{cc} :x^{i\alpha}\partial_{j\beta}:\varphi_{j},& i\not=j\\
    \left(:x^{i\alpha}\partial_{j\beta}:-\sum\limits_{m:m>i}:x^{m\alpha}\partial_{m\beta}:\dfrac{1}{\th_{im}
        \varphi_{im}}\right)\varphi_j,
    & i=j
\end{array}
\right.
\ee
for  fixed $\alpha$ and $\beta$ and then  sum them up over $\alpha=
\beta=1,\ldots,N$.
  First  we note that
  \be\label{O18}
  x^{n\alpha}\scirc \partial_{i\beta}= :x^{n\alpha} \partial_{i\beta}:\qquad\text{and}\qquad x^{i\alpha}\scirc \
  \partial_{n\beta}= :x^{i\alpha}\partial_{n\beta}:
  \ee
  for any $i=1,\ldots ,n$. This is because $x^{n\alpha}$ and $\partial_{n\beta}$ are respectively lowest and highest weight vectors with respect
  to the adjoint action of the diagonal $\Ug$:
  $$ [x^{n\alpha},\e{j}{i}]=0\ ,\  [\e{i}{j},\partial_{n\beta}]=0\qquad\text{for any}\quad i<j,$$
  so that
  $x^{n\alpha}\P\partial_{i\beta}=:x^{n\alpha}\partial_{i\beta}:$ and $x^{i\alpha}\P\partial_{n\beta}=:x^{i\alpha}\partial_{n\beta}:$
  in the double coset space 
  $\Diff(n,N)$.

 Next we apply Zhelobenko operator $\q_{n-1}$ to both sides of equality
$$x^{n\alpha}\scirc \partial_{j\beta}= :x^{n\alpha}\partial_{j\beta}:\,,\qquad j\not=n-1,n.$$
 We have, using \rf{O7} and homomorphism property of Zhelobenko operators,
$$\q_{n-1}(x^{n\alpha}\scirc \partial_{j\beta})=\q_{n-1}(x^{n\alpha)})\scirc\q_{n-1}(\partial_{j\beta})=
x^{n-1,\alpha}\scirc\partial_{j\beta}, \qquad j\not=n-1,n.$$
On the other hand, we can apply $\q_{n-2}$ to $:x^{n-1,\alpha}\partial_{n\beta}:$ as to  elements of the adjoint
representation of $\gl$, see \cite[eq. (4.5)]{KO2}
$$\q_{n-1}(:x^{n\alpha} \partial_{j\beta})=:x^{n-1,\alpha} \partial_{j,\beta}:\,.$$
 This implies, due to \rf{O18}, the equality
$$x^{n-1,\alpha}\scirc \partial_{j\beta}= : x^{n-1,\alpha}\partial_{j\beta}:\qquad \text{for}\quad j\not=n-1.$$
 Proceeding further with application of other Zhelobenko automorphisms $\q_i$ we obtain similarly
\be\label{O19}
x^{i\alpha}\scirc \partial_{j\beta}=  :x^{i\alpha} \partial_{j\beta}:\qquad \text{for}\quad i\not=j\,.
\ee
Due to \rf{O36} this is equivalent to the first line of \rf{O34}. For the derivation of the rest of the relations \rf{O34} we employ 
the action of the Zhelobenko operators on the elements $:x^{i\alpha}\partial_{i\beta}:$ (no sum in $i$). Note that the sum
$x^{i\alpha}\partial_{i\beta}+x^{i+1,\alpha}\partial_{i+1,\beta}$ is invariant with respect to the adjoint action of the
$\bf{sl}_2$-subalgebra generated by $\e{i}{i+1}$ and $\e{i+1}{i}$, while the difference
$x^{i,\alpha}\partial_{i\beta}-x^{i+1,\alpha}\partial_{i+1,\beta}$ spans 
the zero weight subspace of the three-dimensional
representation of this $\bf{sl}_2$-subalgebra, with the highest weight vector $x^{i\alpha}\partial_{i+1,\beta}$.
This implies the relations
\be\label{O20}
\begin{split}
 \q_i\left(:x^{i\alpha}\partial_{i\beta}:\right)&=-:x^{i\alpha}\partial_{i\beta}:\frac{1}{\th_{i,i+1}-1}+
:x^{i+1,\alpha}\partial_{i+1,\beta}:\frac{\th_{i,i+1}}{\th_{i,i+1}-1},\\[-.5em]
 \q_i\left(:x^{i+1,\alpha}\partial_{i+1,\beta}:\right)&=:x^{i\alpha}\partial_{i\beta}:\frac{\th_{i,i+1}}{\th_{i,i+1}-1}-
:x^{i+1,\alpha}\partial_{i+1,\beta}:\frac{1}{\th_{i,i+1}-1}.
\end{split}
\ee
Besides,
\be\label{O21}
\q_i\left(:x^{j\alpha}\partial_{j\beta}:\right)=:x^{j\alpha}\partial_{j\beta}:\qquad\text{for}\quad j\not=i,i+1.
\ee
On the other hand, by \rf{O7} we have
\be\label{O22}
\begin{split}
&\q_i\left(x^{i\alpha}\scirc\partial_{i\beta}\right)=\q_i\left(x^{i\alpha}\right)\scirc\q_i\left(\partial_{i\beta}\right)=
x^{i+1,\alpha}\frac{\th_{i,i+1}}{\th_{i,i+1}-1}\scirc \partial_{i+1,\beta}=
x^{i+1,\alpha}\scirc \partial_{i+1,\beta}\frac{\th_{i,i+1}+1}{\th_{i,i+1}},\\[-.5em]
&\q_i\left(x^{i+1,\alpha}\scirc\partial_{i+1,\beta}\right)=\q_i\left(x^{i+1,\alpha}\right)\scirc\q_{i}
\left(\partial_{i+1,\beta}\right)=
x^{i}\scirc \partial_{i}\,\frac{\th_{i,i+1}}{\th_{i,i+1}-1}.
\end{split}
\ee
We apply $\q_{n-1}$ to both sides of the equality $x^n\scirc \partial_n= :x^n \partial_n:$. Using \rf{O20}
and \rf{O22} we get
$$x^{n-1,\alpha}\scirc\partial_{n-1,\beta}\frac{\th_{n-1,n}}{\th_{n-1,n}-1} =
:x^{n-1,\alpha}\partial_{n-1,\beta}:\frac{\th_{n-1,n}}{\th_{n-1,n}-1}-:x^{n\alpha}\partial_{n\beta}:\frac{1}{\th_{n-1,n}-1},$$
which implies the equality
\be\label{O23}
x^{n-1,\alpha}\scirc\partial_{n-1,\beta}=:x^{n-1,\alpha}\partial_{n-1,\beta}:-:x^{n\alpha}\partial_{n\beta}:\frac{1}{\th_{n-1,n}}.
\ee
Applying $\q_{n-2}$ to \rf{O23} we get
\begin{equation*}
 \begin{split}
&x^{n-2,\alpha}\scirc\partial_{n-2,\beta}\frac{\th_{n-2,n-1}}{\th_{n-2,n-1}-1} =
:x^{n-2,\alpha}\partial_{n-2,\beta}:\frac{\th_{n-2,n-1}}{\th_{n-2,n-1}-1}\\
&\ \ \ \ \ \ -:x^{n-1,\alpha}\partial_{n-1,\beta}:
\frac{1}{\th_{n-2,n-1}-1}-:x^{n\alpha}\partial_{n\beta}:\frac{1}{\th_{n-2,n}},
\end{split}
\end{equation*}
which gives
$$x^{n-2,\alpha}\scirc\partial_{n-2,\beta} =
:x^{n-2,\alpha}\partial_{n-2,\beta}-x^{n-1,\alpha}\partial_{n-1,\beta}
\frac{1}{\th_{n-2,n-1}}-x^{n\alpha}\partial_{n\beta}\frac{\th_{n-2,n-1}-1}{\th_{n-2,n-1}\th_{n-2,n}}:.$$
Proceeding further we obtain for any $i\leq n$ the relation
\be\label{O24}
x^{i\alpha}\scirc\partial_{i\beta} =
:x^{i\alpha}\partial_{i\beta}:-\sum\limits_{m>i}\dfrac{1}{\th_{im}
\varphi_{im}}:x^{m\alpha}\partial_{m\beta}:
\ee
This is precisely the second line of \rf{O34}.\hfill{$\square$}

\vskip .2cm\noindent
{\it Note.} Due to the realization $\LL^i_j\mapsto x^{i\alpha}\scirc \der_{j\alpha}$, the action of the automorphisms $\q_i$ on the generators $\LL^i_j$ can be directly read off
the formulas (\ref{O7a}). In particular, the action on the diagonal generators is standard, $\q_i (\LL^j_j)=\LL^{\sigma_i(j)}_{\sigma_i(j)}$.

\subsection{Central elements}\label{section4.2}

 \paragraph{1.} Let $\Lambda$ be a $n\times n$ matrix with noncommutative entries belonging to some $\Uh$-bimodule, such that the weight of
$\Lambda^i_j$ equals $\ve_i-\ve_j$, that is, $\th_k\Lambda^i_j=\Lambda^i_j(\th_k+\delta^i_k-\delta^j_k)$. Assume that
 $\Lambda$ verifies the reflection equation
 $$\RR_{12}\Lambda_1\RR_{12}\Lambda_1-
 \Lambda_1\RR_{12}\Lambda_1\RR_{12}=\RR_{12}\Lambda_1-
\Lambda_1\RR_{12}\ .$$

\begin{prop}\label{centerth} For any nonnegative integer $N$ the elements $\text{tr}\,(\Lambda^N\QQ^-)$
commute with $\Lambda^i_j$ for all $i$ and $j$.
 \end{prop}
This immediately implies
\begin{cor}\label{corcenterth} For any nonnegative integer $N$ the elements $\text{tr}\,({\L{}{}}^N\QQ^-)$
are central in the algebra $\DR$.
 \end{cor}
\noindent{\it Proof of proposition \ref{centerth}.}
The defining relation
 $$(\RR_{12}\Lambda_1\RR_{12}-\RR_{12})\Lambda_1=\Lambda_1(\RR_{12}\Lambda_1\RR_{12}-\RR_{12})$$ implies that
$$(\RR_{12}\Lambda_1\RR_{12}-\RR_{12})\Lambda_1^N=\Lambda_1^N(\RR_{12}\Lambda_1\RR_{12}-\RR_{12})$$ for any non-negative integer $N$, or, using
(\ref{Q5do}),
\be\label{force1}\Lambda_1\RR_{12}\Lambda_1^N\RR_{12}-\Lambda_1^N\RR_{12}=\RR_{12}\Lambda_1^N\RR_{12}\Lambda_1-
\RR_{12}\Lambda_1^N\ .\ee

Let $\text{Tr}_{(\h)2}$ be the linear map from the space of tensors $\Xi^{i_1i_2}_{j_1j_2}$ to the space of tensors $\Upsilon^{i_1}_{j_1}$ 
such that $\left(\text{Tr}_{(\h)2} (\Xi)\right)^{i_1}_{j_1} =
(\QQ^-)^u_v[-\epsilon_{i_1}]\Xi^{i_1v}_{j_1u}$. This is the $\h$-analogue of the $\R$-matrix trace in the second space. 
We calculate $\text{Tr}_{(\h)2}$ of each term of \rf{force1}.

\noindent$\bullet$ The image of the expression $\Lambda_1\RR_{12}\Lambda_1^N\RR_{12}$ under the map $\text{Tr}_{(\h)2}$ is 
$$\sum_{u,a,b,c,d} \QQ^-_u[-\epsilon_{i_1}]\,\Lambda^{i_1}_{c}\,(\Lambda^N)^a_b\,\RR^{cu}_{ad}[\epsilon_a-\epsilon_b]\,\RR^{bd}_{j_1u}$$
$$=\sum_{u,a,b,c,d} \Lambda^{i_1}_c\,(\Lambda^N)^a_b\,\RR^{cu}_{ad}[\epsilon_a-\epsilon_b]\,\RR^{bd}_{j_1u}\,\QQ^-_u[\epsilon_a-\epsilon_b-\epsilon_c]\ $$
which we rewrite, using (\ref{skew}) and (\ref{skdop1}), as
$$=\sum_{a,b,c}\Lambda^{i_1}_c\,(\Lambda^N)^a_b\,\QQ^-_a[\epsilon_a-\epsilon_b]\,\delta_a^b\,\delta^c_{j_1}=
\Lambda^{i_1}_{j_1}\,\text{Tr}\,(\Lambda^N\QQ^-)\ .$$
$\bullet$ The image of   
$\RR_{12}\Lambda^N_1\RR_{12}\Lambda_1$ is
$$\sum_{u,a,b,d,f} (\Lambda^N)^a_b\,\QQ^-_u[-\epsilon_{i_1}+\epsilon_a-\epsilon_b]\,
\RR^{i_1u}_{ad}[\epsilon_a-\epsilon_b]\,
\RR^{bd}_{fu}\,\Lambda^f_{j_1}\ ,$$
which, using again (\ref{skew}) and (\ref{skdop1}), equals
$$=\sum_{a,b,f} (\Lambda^N)^a_b\,Q^-_a \,\delta_a^b\,\delta^{i_1}_f\,\Lambda^f_{j_1}=\text{Tr}\,(\Lambda^N\QQ^-)\,
\Lambda^{i_1}_{j_1}\ .$$
$\bullet$ The image of $\RR_{12}\Lambda^N_1$ is
$$\sum_{u,a} \QQ^-_u[-\epsilon_{i_1}]\RR^{i_1u}_{au}(\Lambda^N)^a_{j_1}\ .$$
Thus, by (\ref{zadsled}), $\text{Tr}_{(\h)2}(\RR_{12}\Lambda^N_1)=\Lambda^N_1$. 

\noindent$\bullet$ The image of $\Lambda^N_1\RR_{12}$ is
$$\sum_{u,a} \QQ^-_u[-\epsilon_{i_1}]\,(\Lambda^N)^{i_1}_a\,\RR^{au}_{j_1u}=(\Lambda^N)^{i_1}_a\,
\QQ^-_u[-\epsilon_a]\,\RR^{au}_{j_1u}$$
and we again obtain $\Lambda^N_1$.

\noindent Combining these calculations we find that  the application of the map $\text{Tr}_{(\h)2}$ 
to the relation (\ref{force1}) gives $\Lambda\text{Tr}\,(\Lambda^N\QQ^-)=
\text{Tr}\,(\Lambda^N\QQ^-)\Lambda$ as stated.\hfill $\square$

\vskip .2cm
{\it Notes.} 1. $\text{Tr}\,(\LL^N\QQ^-)=\text{Tr}\,(\QQ^-\LL^N)$ since the diagonal elements of $\LL^N$ have weight 0.

\vskip .1cm\noindent 
$\phantom{\it Notes.}$ 2. The reflection equation \rf{O35b} admits shifts $\LL^i_j\to \LL^i_j+\text{const}\cdot\delta^i_j$.

\vskip .1cm\noindent 
$\phantom{\it Notes.}$ 3. $\text{Tr}\,\QQ^+=\text{Tr}\,\QQ^-=n$. Indeed, with the explicit form (\ref{explpsi}) of the tensor $\PPsi$, the relation (\ref{Q2}) is
$$\sum_j\frac{1}{1+\th_{ij}}\QQ^-_j=1\ \text{for any $i$}\ .$$ Write this relation for ${\mathbf {gl}}_{n+1}$, with indices in the range $\{0,1,\dots,n\}$, for $i=0$:
$$\prod_{l=1}^n\frac{\th_{0l}-1}{\th_{0l}}+\sum_{j=1}^n\frac{1}{\th_{0l}}\QQ^-_j=1$$ (here $\QQ^-_j$ corresponds to $\gln$). Decomposing into a power series in
$\th_0^{-1}$ and comparing coefficients at $\frac{1}{\th_0}$ we find $\text{Tr}\,\QQ^-=n$. Since $Q^-_a\vert_{\th\mapsto -\th}=Q^+_a$,
we have $\text{Tr}\,\QQ^+=n$ as well. 

\paragraph{2.} The images, in the reduction algebra, of $\e{i}{j}^{(2)}$ satisfy the same relations as the images of $\e{i}{j}^{(1)}$. 
We have therefore another set of generators of $\DR$
\be\label{Os3bis}
\Lprime^{i}_{j}=\left\{
\begin{array}{cc} \sprime{i}{j}\varphi_{j},& i\not=j\\
    \left(\sprime{i}{j}-\sum\limits_{m:m>i}\sprime{m}{m}\dfrac{1}{\th_{im}
        \varphi_{im}}\right)\varphi_j,
    & i=j
\end{array}
\right.
\ee
which satisfy the same algebra
$$\RR_{12}{\teL}'_1\RR_{12} {\teL}'_1-
 {\teL}'_1\RR_{12}{\teL}'_1\RR_{12}=\RR_{12}{\teL}'_1-
 {\teL}'_1\RR_{12}\ .$$
Using (\ref{relsij}), one can check that the matrices $\teL$ and $\teL'$ are related by
$$\teL'+\teL=\operatorname{H}\ .$$
Here $\operatorname{H}$ is the operator $\C^n\to\Uh\otimes\C^n$ with matrix coefficients
$$\operatorname{H}^i_j:=(\th_j+n)\delta^i_j\ .$$
By proposition \ref{centerth}, for any nonnegative integer $N$ the elements $\text{Tr}\,(\teL'^N\QQ^-)$
are central in the algebra $\DR$.

Substituting into the reflection equation the expression for $\teL'$ in terms of $\teL$ we find
\be\label{Q25}\RR_{12}\operatorname{H}_1\RR_{12}\operatorname{H}_1-
\operatorname{H}_1\RR_{12}\operatorname{H}_1\RR_{12}=\RR_{12}\operatorname{H}_1-
\operatorname{H}_1\RR_{12}\ ,\ee
$$\RR_{12}\teL_1\RR_{12}\operatorname{H}_1+\RR_{12}\operatorname{H}_1\RR_{12}\teL_1-
\teL_1\RR_{12}\operatorname{H}_1\RR_{12}-\operatorname{H}_1\RR_{12}\teL_1\RR_{12}=2\RR_{12}\teL_1-
2\teL_1\RR_{12}\ .$$
Note that the latter equality (one can check it directly) holds for any matrix $\teL$ whose matrix element $\teL^i_j$ has the weight $\varepsilon_i-\varepsilon_j$.

\vskip .2cm
Presumably the center is generated by the elements $\text{Tr}\,(\LL^N\QQ^-)$ and $\text{Tr}\,(\teL'^N\QQ^-)$.

\vskip .2cm
\noindent
{\it Notes}. 1. The center of the algebra of $\h$-deformed differential operators is non-trivial. It is described in \cite{HO}.

\vskip .1cm\noindent
$\phantom{{\it Notes}.}$ 2.
The relation \rf{Q25} shows that the assignment $\LL\to\operatorname{H}$
is a realization of the reflection equation algebra \rf{O35b}.

\subsection{Braided bialgebra structure}\label{section4.3}
Consider the $\Uh$-bimodule $\DR\otimes_{\Uh}\DR$.
The elements $\M^i_j$ generate the first copy of $\DR$ and
$\tilde{\M}^i_j$ generate the second copy. These elements satisfy the relations
\begin{equation}\label{bist1}\begin{array}{ll}
&\RR_{12}\M_1\RR_{12}{\M}_1-\M_1\RR_{12}\M_1\RR_{12}=\RR_{12}\M_1-\M_1\R_{12},\\[.3em]
&\RR_{12}\tilde{\M}_1\RR_{12}\tilde{{\M}}_1-\tilde{\M}_1\RR_{12}\tilde{\M}_1\RR_{12}=
\RR_{12}\tilde{\M}_1-\tilde{\M}_1\R_{12}.
\end{array}\end{equation}
We impose the commutation relation
\begin{equation}\label{bist2}\RR_{12}\M_1\RR_{12}\tilde{\M}_1=\tilde{\M}_1\RR_{12}\M_1\RR_{12}\ .\end{equation}
By virtue of the dynamical Yang--Baxter equation, with this setting,
the $\Uh$-bimodule $\DR\otimes_{\Uh}\DR$ becomes the associative algebra, which we denote
by $\DR\circledcirc\DR $.

In other words, $\DR\circledcirc\DR $ is the associative algebra over $\Uh$, generated
by elements $\M^i_j$ and $\tilde{\M}^i_j$ of weight $\ve_i-\ve_j$ subject to the defining 
relations (\ref{bist1})--(\ref{bist2}). It is isomorphic to $\DR\otimes_{\Uh}\DR$ as a $\Uh$-bimodule.
\begin{lem} The matrix $\M+\tilde{\M}$ satisfies the reflection equation \rf{O35b}. 
\end{lem}
\noindent{\it Proof.} Straightforward.\hfill $\square$
\begin{cor} The map \begin{equation}\label{bist3}\LL\mapsto \M+\tilde{\M}\end{equation} is a homomorphism $\DR\to \DR\circledcirc\DR$
of algebras.
\end{cor}
In a similar fashion we define the product $\circledcirc$ 
of three and more copies of $\DR$. For instance, $\DR\circledcirc\DR\circledcirc\DR$ is generated by 
$\M^i_j$, $\tilde{\M}^i_j$ and $\tilde{\tilde{\M}}^i_j$, with the defining relations (\ref{bist1})--(\ref{bist2}) and, in addition,
$$\begin{array}{c}\RR_{12}\tilde{\tilde{\M}}_1\RR_{12}\tilde{{\tilde{\M}}}_1-\tilde{\tilde{\M}}_1\RR_{12}\tilde{\tilde{\M}}_1\RR_{12}=
\RR_{12}\tilde{\tilde{\M}}_1-\tilde{\tilde{\M}}_1\R_{12},\\
\RR_{12}\M_1\RR_{12}\tilde{\tilde{\M}}_1=\tilde{\tilde{\M}}_1\RR_{12}\M_1\RR_{12}\ \ ,\ \ 
\RR_{12}\tilde{\M}_1\RR_{12}\tilde{\tilde{\M}}_1=\tilde{\tilde{\M}}_1\RR_{12}\tilde{\M}_1\RR_{12}.\end{array}$$
The coproduct  (\ref{bist3}) is clearly coassociative. We have therefore defined the braided bialgebra structure on $\DR$.

\vskip .2cm
The coproduct (\ref{bist3}) has a natural interpretation in terms of differential operators. Divide the interval $\{1,2,\dots,N\}$ into two subintervals, 
$\{1,2,\dots,\nu\}$ and $\{\nu+1,2,\dots,N\}$ for an arbitrary $\nu$, $1\leq \nu <N$. Then $M^i_j\mapsto \sum_{\alpha=1}^{\nu}x^{i\alpha}\scirc\der_{j\alpha}$ 
and $\tilde{M}^i_j\mapsto \sum_{\alpha=\nu+1}^{N}x^{i\alpha}\scirc\der_{j\alpha}$ is a realization of the algebra $\DR\circledcirc\DR$ while 
$L^i_j\mapsto \sum_{\alpha=1}^{N}x^{i\alpha}\scirc\der_{j\alpha}$ is a realization of the algebra $\DR$.
\subsubsection*{Appendix. Basic relations for $n=2$}
 Denote $\th=\th_{12}=h_1-h_2+1$.
 The defining relations  between different copies of generators in the algebra $\Diff(2,2)$, see \rf{O31}, 
   look as follows:
  \begin{equation*}
  \begin{aligned}
  	x^{1}\scirc {x'}^{2}&=\frac{1}{\th_{}}{x'}^{1}\scirc {x}^{2}+
  	\frac{\th_{}^2-1}{\th_{}^2}{x'}^{2}\scirc x^{1},\qquad&
  	\der_{1}\scirc {\der'}_{2}&=-\frac{1}{\th_{}}{\der'}_{1}\scirc {\der}_{2}+
  	\frac{\th_{}^2-1}{\th_{}^2}{\der'}_{2}\scirc \der_{1},\\[-.6em]
  	x^{2}\scirc {x'}^{1}&={x'}^{1}\scirc x^{2}-\frac{1}{\th_{}}{x'}^{2}\scirc x^{1}
  	,\, &
  		{\der}_{2}\scirc {\der'}_{1}&={\der'}_{1}\scirc \der_{2}+\frac{1}{\th_{}}{\der'}_{2}\scirc \der_{1}
  		, \\[-.6em]
  		x^i\scirc {x'}^i&={x'}^i\scirc {x}^i,\qquad&
  			\der_i\scirc{\der'}_i&={\der'}_i\scirc\der_i,\qquad\qquad i=1,2,
  \end{aligned}
  \end{equation*}
 \begin{equation*}
  	\begin{aligned}  
  		x^1\scirc {\der'}_2&={\der'}_2\scirc x^1,\qquad&
  		x^2\scirc {\der'}_1&=\frac{\th_{}(\th_{}+2)}{(\th_{}+1)^2}\,{\der'}_1\scirc x^2,
  		\\[-.6em]
  		x^1\scirc {\der'}_1&= {\der'}_1\scirc x^1+\frac{1}{1-\th_{}}\,{\der'}_2\scirc x^2-1,\qquad&
  		x^2\scirc {\der'}_2&= \frac{1}{1+\th_{}}\,{\der'}_1\scirc x^1
  		+{\der'}_2\scirc x^2-1.
  	\end{aligned}
\end{equation*} 
Here the elements $\{x^1,x^2,{\der}_1,{\der}_2\}$ belong to the first copy and the elements $\{{x'}^1,{x'}^2,{\der'}_1,{\der'}_2\}$ belong to the second copy. 

The ordering form of the relations \rf{O35b} is
\begin{align*}\LL^1_1\LL^1_2&=\frac{\th-3}{\th-2}\LL^1_2\LL^1_1+\frac{1}{\th-2}\LL^1_2\LL^2_2+\LL^1_2\ ,\\[-.2em]
\LL^2_2\LL^1_2&=\frac{\th-3}{(\th-2)(\th+1)}\LL^1_2\LL^1_1+\frac{(\th-1)^2}{(\th-2)(\th+1)}\LL^1_2\LL^2_2-\frac{\th-1}{\th+1}\LL^1_2\ ,\\[-.2em]
\LL^1_1\LL^2_1&=\frac{(\th+1)^2}{(\th-1)(\th+2)}\LL^2_1\LL^1_1-\frac{\th+3}{(\th-1)(\th+2)}\LL^2_1\LL^2_2-\frac{\th+1}{\th-1}\LL^2_1\ ,\\[-.2em]
\LL^2_2\LL^2_1&=-\frac{1}{\th+2}\LL^2_1\LL^1_1+\frac{\th+3}{\th+2}\LL^2_1\LL^2_2+\LL^2_1\ ,\\[-.2em]
\LL^1_1\LL^2_2&=\LL^2_2\LL^1_1\ ,\\[-.2em]
\LL^1_2\LL^2_1&=\LL^2_1\LL^1_2-\frac{1}{\th}(\LL^1_1-\LL^2_2)^2+\LL^1_1-L^2_2\ .
\end{align*}
The central elements of Corollary 4.4 have the form
$$\frac{\th-1}{\th}(\LL^N)^1_1+\frac{\th+1}{\th}(\LL^N)^2_2\ .$$ 

\vskip .2cm
\noindent{\bf Acknowledgments.} The work of S. K. was supported within the framework 
of a subsidy granted to the HSE by the Government of the
Russian Federation for the implementation of the Global Competitiveness Program, by the RFBR grant 14-01-00547 and joint RFBR-Ukrainian Academy of Sciences grant 14-01-90405-Ukr.

\end{document}